\documentclass[12pt]{article}

\usepackage{amsmath}
\usepackage{amssymb}
\usepackage{amsthm}
\usepackage{amscd}
\usepackage{amsfonts}
\usepackage{dsfont}
\usepackage{cases}
\usepackage[applemac]{inputenc}
\usepackage{enumerate}
\usepackage[all]{xy}
\usepackage{esint}
\usepackage{graphicx}
\usepackage{url}
\theoremstyle{plain}
\newtheorem{theor}{Theorem}[section]
\newtheorem{Lem}[theor]{Lemma}

\theoremstyle{definition}

\newtheorem{rem}[theor]{Remark}

\newcommand{\N}{\mathbb N}

\newcommand{\R}{\mathbb R}

\newcommand{\E}{\mathrm E}
\newcommand{\Ess}{\cal S}
\newcommand{\T}{{\cal T}}

\newcommand{\D}{\Delta}
\newcommand{\G}{{\cal G}}
\newcommand{\F}{{\cal F}}
\usepackage[authoryear]{natbib}
\usepackage{url}
\date{}
\begin{document}
\title{On the $1/e$-strategy for the best-choice problem under no information.}
\author{F. Thomas Bruss\\Universit\'e Libre de Bruxelles}
\maketitle

\noindent{\bf Abstract~} The main purpose of this paper is to correct an error in the previously submitted version [*] := arXiv:2004.13749v1. [*] had been already accepted for  publication in a scientific journal, but withdrawn by the author after the discovery of the error. For the withdrawal from arxiv we follow their preference to maintain what remains of interest.
The background of the open problem, and the brief survey which comes with it, stay relevant. These keep their place in the present corrected version. The same is true for two new modified  odds-theorems proved in [*] since they are applicable for several different stopping problems.  Then, and in particular,
 we show where exactly the error occurred in [*],  why it invalidates its main theorem and title, and what the conclusions are. The final discussion of optimal strategies {\it without value} in Section 4 is believed to be of general independent interest.


\bigskip\noindent{\bf MSC 2010 Subject Code}: 60G40
\section{Background of the problem}
At the evening of Professor  Larry Shepp's talk ``Reflecting Brownian Motion" at Cornell University on July 11, 1983 (13th Conference on Stochastic Processes and Applications), Professor Shepp and the author ran into each other in front of the Ezra Cornell statue. I was  honoured to meet him in person, and Larry replied ``What are you working on?" And so Prof. Shepp was the very first person with whom I could discuss the {\it $1/e$-law of best choice} resulting from the {\it Unified Approach} (B. (1984)) which had been accepted for publication shortly before. 
I was glad to see his true interest in the $1/e$-law. As many of us  know, when Larry was interested in a problem, then he was deeply interested.\smallskip

This article deals with an open question concerning the optimality of the so-called $1/e$-strategy for the problem of best choice under no information on the number $N$ of options. I drew again attention to this open question in my own talk ``The $e^{-1}$-law in best choice problems" at Cornell on July 14, 1983, and re-discussed it with Larry at several later occasions. 
An earlier related question appears already on page 885 of B. (1984) where the author conjectured that, for a {\it two-person game}, the $e^{-1}$-strategy is optimal for the decision maker who has to select. As far as the author is aware, the last written reference to the precise open question discussed with Prof. Shepp is in B. and Yor (2012).
\smallskip

\section {The Unified Approach}
We begin with a  review of the {\it Unified Approach}-model and previously known results.
\begin{quote}{\bf Unified Approach}: Suppose $N>0$ points are i.i.d. with a continuous distribution function $F$ on some interval $[0,T].$
Points are marked with qualities which are supposed to be uniquely rankable from $1$ (best) to $N$ (worst), and all rank arrival orders are supposed to be equally likely. The goal is  to maximize the probability of stopping online and without recall on rank $1.$ (B. (1984)) \end{quote}
\noindent This model  was suggested for the best choice problem (secretary problem) for an unknown number $N$ of candidates. Recall that, by R\'enyi's theorem of relative ranks (R\'enyi (1962)), the $k$th candidate has relative rank $j$ with probability $1/k$ for all  $1\le j\le k$ whenever all rank arrival orders are equally likely.

Previous models for unknown $N$ had shown that the price for not knowing $N$ can be high. The influential paper by 
Presman and Sonin (1972)  which modelled the unknown $N$ via the hypothesis of a known distribution $\{P(N=n)\},$ displayed the intricacies arising by the possible appearance of so-called {\it stopping islands.} Moreover, Abdel-Hamid et al. (1982) showed that the $N$-unknown problem may have several solutions, and, much worse, that for any $\epsilon>0$ there exists a sufficiently unfavorable distribution $\{P(N=n)\}_{n=1,2, \cdots}$ to reduce the optimal success probability to a value smaller than   $\epsilon.$ In other words, if $N$ is modelled via $\{P(N=n)\},$ optimality may mean in some cases almost nothing. This contrasts with the well-known lower bound $1/e$ which holds in the classical model for known $N=n\ge 1.$
These discouraging facts  for unknown $N$ instigated efforts to find more tractable models, as e.g.\,the model of Stewart (1981), or the one of Cowan and Zabzcyk (1978) and its generalisation studied in B. (1987), and also others.

\smallskip The unified approach of B. (1984) was different.  The idea behind it was that it is typically easier to estimate - and this is where the  time distribution $F$ comes in - {\it when} options  are more likely to arrive conditional on knowing that they do arrive than making hypotheses about the
distribution of its {\it number.} No assumption at all was made about the distribution of $N.$ (The same approach was later taken by B. and Samuels (1987) for more general payoffs for different ranks.)

\smallskip 
The continuous arrival time distribution $F$ is the crucial part with respect to applications.
For our open problem itself the form of $F$ is irrelevant, however. If we transform the unordered i.i.d arrival times of the best, the second best ... , $T_1, T_2, ... $ say, by $T^*_k:=F(T_k)$, then the $T^*_k$ are i.i.d. $U[0,1]$ random variables and, since $F$ is continuous and increasing, the time transformation maintains the arrival order of the different relative ranks. Thus, if we know the optimal strategy for dealing with i.i.d. $U[0,1]$ random arrivals on $[0,1]$, then we know it as well for i.i.d. $F$-distributed
arrival times on the original horizon $[0,T].$ In all what follows we therefore confine our interest to uniformly distributed arrival times in $[0,1]$-time.
\subsection{ Related problems}
\smallskip
A  related problem, to which we will return in Subsection 2.6, is the so-called last-arrival-problem
under no information (l.a.p.)\,studied by B. and Yor (2012).

In this model an unknown number $N$ of points are  i.i.d. $U[0,1]$ random variables, and an observer, inspecting the interval $[0,1]$ sequentially from left to right, wants to maximise the probability of stopping on the very last point. No information about $N$ whatsoever is given.
Only one stop is allowed, and this again without recall on  preceding observations (online). Thus the only difference of the l.a.p. model of B. and Yor (2012) to the Unified Approach model of B. (1984) is that no ranks are  attributed to the observations (points).

\medskip Other related problems, now again with the objective to get rank 1 of uniquely ranked candidates, arise by combining the Unified Approach model and the model of Presman and Sonin (1972) for different types of distributions of $N$. If $(P(N=n))_{n=1,2, \cdots}$ is known, then one is in the setting of a model with a prior. The i.i.d. $U[0,1]$ arrival-times can then be used as an additional means of statistical inference to update the posterior distribution of $N$. Stopping islands, as observed  in the paper of Presman and Sonin (1972), bear over to corresponding islands in continuous time. The optimal strategy may thus become very complicated, and we would typically not like to compute it, but, in principle, it can be computed.

\medskip For the latter class of problems, what would be a good alternative? Moreover, and in particular, what can one do if one has absolutely no information about $N$?

\subsection{The 1/e-law}
The answer given by the unified approach (B.(1984)) was  that, as far as  applications are concerned,  we need not care much.    For ease of reference we recall these results summarised as the {\it $1/e$-law.} Here we follow the meanwhile established tradition to call an observation of relative rank 1 a {\it record value}, or simply {\it record,} and the time when a record appears a {\it record time.} R\'enyi (1962) had called a record an {\it \'el\'ement saillant}.

\medskip
The $1/e$-law says:
\begin{quote}
1. The strategy to wait (in $[0,1]$-time) up to time $1/e \approx 0.3678,$ and then to select the first record (if any from time $1/e$ onward), called the $1/e$-strategy, succeeds for all $N$ with
probability at least $1/e.$

\smallskip
2. There exists no strategy which would be better for all $N.$

\medskip
3. The $1/e$-strategy selects no candidate with precise probability $1/e.$\end{quote}

\noindent  Note also that 1. and 3. imply that a non-best option is selected with probability smaller than $1-2/e\approx 0.2642.$
This multiple role of the number $1/e$ gave rise to the name $1/e$-law, and Table 1 (B.(1984), p. 336) shows how good the lower bound $1/e$ for the success probability actually is.
Taking also into account the minimax optimality stated in 2. we can conclude that the $1/e$-strategy is   a convenient and convincing alternative
for all practical purposes.
See e.g. the comments of Samuels (Math. Reviews: 1985). 

\smallskip

But then, the following question arises:
\smallskip

{\it Is the $1/e$-strategy optimal 
if we have no prior information at all on $N$?}

\medskip
\noindent As mentioned before, if the question is stated like this the answer is No. We have to  return to what is known.

\medskip
What is known?

\medskip
(I)~{\bf Optimal $x$-strategies given $N=n.$}~ First, suppose that $N$ were known, say $N=n$, and that we want to determine the optimal strategy in the class of so-called $x$-strategies, that is to wait  until time $x\in [0,1]$ and then to select, if any, the first record from time $x$ onward. 
It is not difficult to find, conditioned on $\{N=n\}$, the optimal waiting time $x_n$ and its performance in this class of $x$-strategies,  namely (see B. (1984), p. 884, (2)-(7)),
\begin{align}  ~~x_1=0; ~x_n=\arg\left\{0 \le x \le 1: \sum_{k=1}^{n-1}\frac{(1-x)^k}{k} = 1\right\}, n=2,3, \cdots.\end{align} 
Note that the $x_n$-strategy is suboptimal since it does not fully use the knowledge $N=n,$ as it is the case for the optimal strategy for the classical secretary problem for $n$ candidates.

\bigskip

(II)~{\bf Monotonicity results.~} We can derive from (1) that  \begin{align}p_n(x):= P(x{\rm -strategy~succeeds}\big |\,N=n)=\frac{(1-x)^n}{n}+x\,\sum_{k=1}^{n-1}\frac{(1-x)^k}{k},\end{align}
and also that  $p_n(x)\ge p_{n+1}(x)$ for all $x\in [0,1].$ This implies
\begin{align}
\forall x \in [0,1]: p_n(x) \downarrow \lim_{n\to \infty}p_n(x)=- x \log (x).
 \end{align}Moreover, it follows from (2) and (3) that the optimal waiting time $x_n$ and the corresponding optimal win probability $p_n(x_n)$ satisfy, respectively, \begin{align}~~x_n\uparrow \frac{1}{e} ~~{\rm and~~}p(x_n)\downarrow \frac{1}{e}, {~\rm as~}n \to \infty.\end{align} 

\medskip 

(III) {\bf Asymptotic optimality.}~The  $1/e$-strategy is, as $n\to \infty$, asymptotically 
optimal with win probability $1/e.$ This follows from (3) and (4), showing that the limiting performance of the $1/e$-strategy is the same as that of the well-known optimal strategy for  the  classical secretary problem for known $n$ as $n\to \infty,$ namely $1/e.$ Clearly one cannot do better than in the case that one knows $N.$

\bigskip (IV) {\bf Connection with Pascal-processes}
~ Let $(\Pi_t)_{t\ge0}$ be a counting process on $\R^+$ with the distributional prescription that for all $T>0$ and $0< t\le T$
\begin{equation*}P(\Pi_T=n|{\cal F}_t)={n \choose \Pi_t}p(t,T)^{\Pi_t +1}(1-p(t,T))^{n-\Pi_t},
\end{equation*}
where $\Pi_0=0$ and $({\cal F}_t)=\sigma(\{\Pi_u:u\le t\}).$ Then $(\Pi_t)$ is called a Pascal process with parameter function $p(t,T).$ These processes are characterized in B. and Rogers (1991).
Pascal processes have the remarkable property that  if points are marked independently with ranks, then, concentrating on 1-records ($\equiv$ records) in such a process, 
optimality for stopping on the last record cannot depend on the number of points
seen before. This property of stationarity was earlier observed in a weaker form (quasi-stationarity) in B. and Samuels (1990). Both papers thus add to the interest of knowing the answer of the open problem.

\bigskip
{\bf Challenge and Intuition}
\smallskip

The mathematical challenge to have a complete answer for the case of no-information 
remains because the unified approach model was created in order to deal with any $N.$  What attempts were made before, and why?

Looking in (II) of Subsection 2.2 closely at (2), (3) and (4),  the open problem comes up quite naturally. Things become intriguing. For any $N=n$ there is a better strategy since the optimal $x_n$-waiting time strategy turns out strictly better than the $1/e$-strategy. Thus one gets the feeling that if there were a way of collecting information about $N$ sufficiently quickly, then this may be sufficient to prove that the $1/e$-strategy cannot be optimal. 
Viewing to disprove optimality, it seems promising to assume  certain types and amounts of weak information about $N$, still strong enough to imply that the $1/e$-strategy is {\it not} optimal, and then to weaken the information. 

Interestingly, as soon as one lets information about $N$ become weaker and weaker, and finally fade away towards no-information, the $1/e$-waiting time seems to become a miraculous  "fix-point" of optimal thresholds.  According to III, this would surprise us much less if no-information on $N$ implied in any way that $N$ is likely to be large, but of course it does not! To understand this is a challenge.

\smallskip
What about trying to find other types of counterexamples?
\smallskip

The challenge remains. It is not easy to do this without leaving the framework of no-information.  Arguing for example ``If we have no information on $\{P(N=n)\}_{n=1, 2, \cdots},$ then let us for instance suppose that this distribution turns out such and such, and that we have seen a history of points such and such, ..." and then imply that the $1/e$-strategy is {\it not} optimal, would not be correct. 
Proofs by contradiction are only valid within the same logical framework, i.e.\,no-information. Arguments
implying initial information whatsoever on $N$ would not be rigorous. For the same reason, simulations are meaningless as they require parameters to  randomize $N,$ and thus information on $N$ must be inputted.
Looking for counterexamples cannot be expected to help. Knowing this increases the challenge.

\subsection{Ill-posed or well-posed problem?} \label{sec26}

Is the question possibly ill-posed?
This question was asked repeatedly by several peers, and, during certain periods, the author also shared these doubts. 
Indeed, the notion of ``no-information" requires clarification. 

Can one interpret no-information in the sense that all possible values of $N$ are in an {\it unknown} interval $\{1,2, \cdots, n\}$ with no value of $N$ being more likely than others, and then let $n$ tend to infinity?  

No. This is equivalent to the improper Laplace prior for $N$. It is true that this prior is the prime candidate for no-information, and very often used to express the lack of knowledge about a parameter. However, this prior implies that $N$ is likely to be very large, and this \textit{is} information.  After all,  ``no information" on $N$ should mean that  at time $0$ we know really  nothing at all about $N.$

Now, more importantly, can we assure that the no-information hypothesis is a honest hypothesis in the sense that it is contradiction-free? If it is not contradiction-free, then of course we must declare  the open problem ill-posed.

\subsection {Formalising no-information}
When B. and Yor (2012) studied the last-arrival problem (l.a.p.) under the no-information hypothesis, they faced a similar difficulty of knowing whether their problem is well-posed.  These authors found a simple argument to prove that it is {\it impossible to prove} that the no-information hypothesis may imply contradictions! 
Their argument was that, whatever a hypothetical information space $\cal H$ may look like for the unknown parameter or random variable $N$, no-information is bound to be a singleton in that space $\cal H$. 

This definition may first sound like a formalism  to prevent saying something wrong, but there is more to it.  It implies that, as a singleton, the no-information hypothesis cannot lead  to contradictory implications. A singleton in $\cal H$ has, by definition of a singleton,  nothing in common with other points in $\cal H,$ whereas contradicting implications cannot come out of nothing. They would need different sources of information giving rise to (at least two) different implications. 

B. and Yor (2012) therefore concluded that they should, a priori, take a constructive attitude and try to find a solution. And so they did. But then the question was to know whether their solution is the solution of a well-posed problem. Hadamard's criteria  (Hadamard (1902))  were the only criteria Bruss and Yor found for the solution of a well-posed problem, and they agreed with these criteria. This is why they were glad to see that their solution fully satisfied Hadamard's criteria.  According to these  criteria, the solution of Bruss and Yor (2012,  subsection 5.3) is the solution of a well-posed problem.

\smallskip
One part of the approach of B. and Y. (2012),
following next, remains however very helpful for our problem, and this is the notion of a stochastic process with proportional increments.

\subsection{ Proportional increments}
For $N$ i.i.d. $U([0,1])$ arrival points, let   \begin{align*}N_t=  \#~ {\rm arrivals ~up~ to~ time~} t,~ t\in[0,1]. \end{align*} B. and Yor (2012, subsection 1.1 and  pp. 3242-46) showed that the counting process $(N_t)_{0\le t\le 1}$ of incoming points on $[0,1]$ with $N:= N_1$
 can be seen as a history-driven process with, what they called, {\it proportional increments}.  This means
that the process $(N_t)$ must satisfy \begin{align*}\forall \,0<t<1 {\rm ~with~}N_t>0:~~~~~~~~~~~~~~\\E(N_{t+\Delta t}-N_t|{\cal F}_t)=\frac{\Delta t}{t} N_t ~a.s., ~ 0 < \Delta t \le 1-t,\end{align*} where the condition $N_t>0$ is crucial. Such a process $(N_t)$ will be said to have  the property of proportional increments, or in short, having the {\it p.i.-property.}

Conditioned on $N=N_1>0,$ let $T_1$ be the first arrival time. The definition of the p.i.-property implies then that, given $N>0,$ the process $(N_t/t)$ is a martingale on $[T_1,1]$, as shown in B. and Yor (2012, p. 3245). This clearly holds also under the stronger assumption that $(N_t)$ is a Poisson process on $[0,1].$  However, B. and Yor (2012, see  p. 3255)  saw a true benefit in {\it not} imposing that $(N_u)$ be Poisson.

To be complete on this, we should mention that in our open problem we could, from a purely decision-theoretic point of view, suppose right away  that the process $(N_t)_{0\le t\le 1}$ is a Poisson point process with unknown rate. Indeed, this cannot make a difference for decisions because we cannot tell a counting process which leaves a pattern of arrival times of a homogeneous Poisson process with unknown rate from another counting process leaving, in distribution, the same pattern. Doing so would have the advantage to be able to use the same compensator on the whole interval $[0,1].$  However, we will not need the Poisson process assumption.

\subsection {Towards suitable odds-theorems}
Recall that our open problem is different from the l.a.p. of B. and Yor (2012)  since  in the Unified Approach model we would like to stop on the very last record, not on the very last point, and thus our approach must also be different. 

The very first arrival time $T_1$ in the counting process $(N_t)$, is the time when $(N_t)$ makes its first jump, and the processes $(N_t)$ and $(N_t/t)$ have exactly the same jump times.  $T_1$ is also the birth time of the record process $(R_t)$, say, defined by $$R_t=\# ~\rm {records~on~[0,t]}, ~ 0\le t \le 1.$$
Since in our open problem any strategy is equivalent for $N=0$ we may and do suppose that $N>0$ almost surely, and thus $T_1<1$ almost surely. $N$ is unknown at time $0$, but at time $1$ we know that, by definition, $N=N_1\ge 1$ almost surely. 

Since the first arrival is also the first record,  we have $N_{T_1-}=R_{T_1-}= 0$ and $N_{T_1}=R_{T_1}= 1.$  Thus the two processes $(N_t)$ and $(R_t)$ have the same random birth time $T_1.$ On the interval $[T_1,1],$ the process $(N_t)$ has proportional increments, i.e. dependent increments, whereas $(R_t)$ has, as we shall see later, under certain conditions independent non-homogeneous increments. We recall here that, by the no-information hypothesis, we have no access to the posterior distribution $\{P(N=n\,|\,T_1=t_1)\}_{n=1, 2, \cdots}.$

To prepare for these properties of $(R_t),$ the idea is to first concentrate on its increments (after time $T_1$). For this purpose we prove two suitably extended versions the Odds-Theorem of optimal stopping. We should also mention here that Ferguson (2016) gave several interesting extensions of the Odds-Theorem in other directions. Moreover, Matsui and Ano (2016) studied in another extension lower bounds of the optimal  success probability for the case of multiple stops. However, here we will here new extensions which are specifically tailored for our open problem. 
We begin with an extension in discrete time.

\subsection{Odds-Theorems for delayed  stopping}
 Let $n$ be a positive integer,  and let $X_1, X_2, \cdots, X_n$ be independent Bernoulli random variables with success parameters $p_k=P(X_k=1)=1-P(X_k=0), k = 1, 2, \cdots, n.$ Suppose our goal is to maximize the probability of stopping online on the last success, i.e. on the last $X_k=1.$ 
The optimal strategy to achieve this goal is immediate from the Odds-Theorem (Bruss (2000)) which we recall for convenience of reference. Let 
\begin{align} q_k=1-p_k;~ r_k=\frac{p_k}{q_k};~ R(k,n) = \sum_{j=k}^n r_j, ~k=1,2,\cdots,n,\end{align}
and let the integer $s\ge1$ (called threshold index) be defined by\begin{align}
 s =
  \begin{cases}
   1 & ,\text{if~} R(1,n) <1 \\
  \sup\{1\le k\le n: R(k,n) \ge 1\}& , \text{otherwise.} \\
     \end{cases}
\end{align}The  strategy to stop  on the first index $k$ with $k\ge s$ and $X_k=1$ (if such a $k$ exists) maximises the probability of stopping on the very last success (B. 2000). If no such $k$ exists, then it is understood that we have to stop at time $n$ and lose by definition. (For a different payoff function and a different approach see e.g. Grau Ribas (2020).) 

\bigskip\noindent
{\bf Delayed stopping in discrete time}
  
  \smallskip
\noindent Let us now consider the new case that there is a deterministic or a random delay imposed by a random time $W$ with values in $\{1, 2,\cdots,n\}$ in the sense that stopping is not allowed before time $W$. Our objective, as before, is to maximize the probability of stopping on the last success.
 Does it suffice to replace simply the threshold $s$
 defined in (6) by $\tilde s:=\max\{W,s\}$ to obtain an optimal strategy? This seems trivial (and is true) if $W$ is deterministic. 
 
 In general this is not true, of course, even not true if $W$ is a stopping time on $X_1, X_2, \cdots, X_n,$ unless we can guarantee that the knowledge of $W$ has no effect on the {\it laws} of $X_{W+1}, X_{W+2}, \dots$ and their independence. 
 The following is a more tractable formulation.\medskip

\begin{theor}  {\it  Let $X_1, X_2, \cdots, X_n$ be Bernoulli random variables defined on a filtered probability space $(\Omega,{\cal A}, ({\cal A}_k), P)$ where ${\cal A}_k = \sigma(\{X_j: 1\le j \le k\}).$ Suppose there exists a random time $W$ for $X_1, X_2, \cdots, X_n$ on the same probability space  such that  the $X_j$ with $j\ge W$ are independent random variables satisfying \begin{align*}p_j(w):=P(X_j=1|W\le w), ~1 \le w\le j\le n. \end{align*}
Then, putting $r_j(w)=p_j(w)/(1-p_j(w)), $ it is optimal to stop at the random time
\begin{align}\tau:=\inf\left\{k \in \{W, W+1,\cdots,n\}: \{X_k=1\} ~\& ~\sum_{j=k+1}^n r_j(W) \le1\right\},\end{align} with the understanding that we stop at time $n$ and lose by definition, if  $\{...\}=\emptyset.$}\end{theor}

\begin{rem} We note that no (initial) independence hypothesis is assumed for the $X_1, X_2, \cdots$ but only for those $X_j$'s with $j\ge W.$ \end{rem}

\medskip \noindent {\bf Proof~of Theorem 2.1}  Our proof will profit from the proof of the Odds-Theorem (B. (2000)) if we rewrite the threshold index (6) in an equivalent form.

\smallskip\noindent Recall the definition of $R(k,n)$ in (5). If we define, as usual, an empty sum as zero, then  $s$ defined in (6) can be written as \begin{align}s'=\inf \left\{1\le k \le n: R(k+1,n):=\sum_{j=k+1}^n r_j\,\le 1\right\}.\end{align} This is straightforward: If $R(1,n)\le1$ then $R(2,n)\le1$ so that from (8) $s'=1$, and $s=1$, as stated in (6). Otherwise, if $R(1,n)> 1,$ then there exists  a unique $k$ where $R(k+1,n)$ drops for the first time below the value $1$ since $R(k,n)$ decreases in $k$, and $R(n+1,n)=0.$ The first such $k$ is the $s'$ defined in (8).
The definitions  (6) and (8) are thus equivalent. (See also Stirzaker (2015, p. 50))
\smallskip

Let now $p_j(w)$ as defined in Theorem 2.1, and let for $ j= w, w+1, \cdots, n$$$q_j(w) = 1-p_j(w)=P(X_j=0|W\le w).$$It follows from the assumptions concerning $W$  that $X_{w}, X_{w+1}, \cdots , X_n$ are independent random variables with laws only dependent on the event $\{W\le w\}.$ If we think of $w$ as being fixed, then we can and do define $p_j:=p_j(w)$ for all $w\le j\le n$ and use the same notation as before defined in (5) with the corresponding odds $r_j(w)=p_j(w)/q_j(w)=: r_j.$
  Accordingly, we have  for $k\ge w$ the same simple monotonicity property $R(k,n)\ge R(k+1,n).$ 
  
  It is easy to check that this monotonicity property is equivalent to the uni-modality property
proved in B. (2000, p.1386, lines 3-12). The latter implies that the optimal rule is a {\it monotone} rule in the sense that, once it is optimal to stop on a success at index $k,$ then it is also optimal to stop on a success after index $k.$ (For a convenient criterion for a stopping rule in the discrete setting being monotone, see Ferguson (2016, p. 49)). 

\smallskip
 Note that, whatever $W=w\in \{1,2, \cdots, n\}$, the odds $r_j:=r_j(w)$ are deterministic functions of the $p_j:=p_j(w)$,  and so the future odds $(r_j )_{j\ge W+1}$ are also known and will not change. The only restriction we have to keep in mind for the simplified notation is that $k\ge w$ on the set $\{W\le w\}.$ 
 But then the monotonicity property of $R(\cdot,\cdot)$ is also not affected, that is
$$\forall \ell \ge j: R(W+j,n)\le1 \implies R(W+\ell,n)\le1.$$  Since the latter implies the uni-modality property of the resulting win probability on $W\le j\le n$,
 the monotone rule property is again maintained for the optimal rule after the random time $W,$ exactly as in B. (2000). Therefore the optimal strategy is to stop on the first success (if it exists) from time $\tau$ onwards where $\tau$ satisfies \begin{align} \tau \ge W ~{\rm and ~}\sum_{j=\tau+1}^n r_j(W) \le1. \end{align} This is the threshold index $\tau$ of Theorem 2.1, and hence the proof.\qed

\medskip
\begin{rem}
Note that Theorem 2.1 is intuitive. Its applicability, nevertheless, can be delicate. It depends on the $p_j$'s being predictable for {\it all} $j\ge W.$   Often this is not the case. For instance, we may have (conditionally) independent random variables, but, if we collect information about the $p_j$ from observations then the distributions of the future values of $X_{j+1}, X_{j+2}, \cdots$ typically depend  on $X_k,~ {1\le k\le j},$ on which the event $\{W=j\}$ may be allowed to depend!
(For our purpose of settling the open question  the implications of Theorem 2.1 will turn out to be strong, however.)
\end{rem}

\begin{rem} (Side-remark). Given that (8) is a one-line definition whereas (6) needs two lines, some readers ask why B. (2000) used definition (6). The answer is that it is (6) which points to the odds-algorithm 
(subsection 2.1, p.1386) which works backwards until the stopping time $s$ with $r_n, r_{n-1}, \cdots$ to give both optimal strategy and value at the same time. No other algorithm can be quicker since it computes exactly those  $r_j$ which produce both answers.
If we used instead the odds beginning with $r_1, r_2, \cdots$ and (8) we would first need $R(1,n)$, implying in general redundant calculations. For the preceding theorem, however, we clearly needed (8). \end{rem} 

\medskip
\noindent{\bf Delayed stopping in continuous time}

\medskip 
\noindent We now state and prove a continuous-time analogue  of the Theorem 2.1 which plays an important role in the proof  of the open conjecture. We  state and prove it in a slightly more general form than what we need  for the conjecture, because it may be also of interest for other problems of optimal stopping.
\begin{theor} Suppose $(C_t)$ is a  counting process on $[0,1]$ for which there exists a random time $\cal T$ such that the confined process $(C_t)_{{\cal T}<t\le 1}$ has independent increments according to a predictable (non-random) intensity measure $\eta(t)_{{\cal T}< t\le1}.$
We suppose that $\eta(t)$ is Riemann integrable on $[0,1]$ with $\E(C_1)<\infty$. Then the optimal strategy to stop on the {\it last} jump-time of $(C_t)$ is to select, if it exists, the first arrival time $\tau \ge \cal T$  with $\tau$ satisfying
\begin{align} \E(C_1-C_\tau)= \int_{\tau}^1 \eta(u) \le 1.\end{align}
\end{theor}

\medskip\noindent 

We note that when the process $(C_t)$ is a Poisson process on $[0, 1]$ the conditions of Theorem 2.5 are clearly satisfied everywhere on $[0,1].$  This special case has been studied already in subsection 4.1 of B. (2000). 

\bigskip\noindent{\bf Proof of Theorem 2.5}

\smallskip\noindent
Consider a partition $\{u_0< u_1< \cdots <u_m\},~ m \in \{1,2, \cdots\},$ of the random sub-interval $[{\cal T},1]\subseteq [0,1]$ with $u_0={\cal T}$ and $u_m=1$. Let  the index $j$ be thought of as depending on $m,$ thus $j:=j(m)$ and $u_j:=u_{j(m)}.$ Put
\begin{align}p_j:=p_{j(m)} = \int_{u_{j -1}}^{u_j}\eta(u) du, ~ j=1, 2, \cdots, m, \end{align} where $[u_{j-1}, u_j[$ is by definition the $j$th sub-interval of the partition, $j=1, 2, \cdots, m.$ It follows that $p_j$ is the expected number of points of the process $(C_u)$
in the $j$th sub-interval, and thus by additivity from (11)\begin{align} \sum_{j=1}^m p_j= \int_{\T}^1 \eta(u)du = \E(C_1-C_\T)\le \E(C_1) <\infty.\end{align} Since all $p_j$ in (11) are non-negative, and $\E(C_1)$ is finite, we can interpret them all as probabilities of certain events as soon as we choose sufficiently fine partitions to have the the $p_j$ less or equal to $1.$ This is always possible since, as we see in (11), $p_j \to 0$ as $\Delta_j=u_j-u_{j-1}\to 0.$ For the following it is understood that we only speak of such sufficiently fine partitions.
Since the counting process $(C_u)$ has independent increments, this allows us at the same time to see the $p_j$ as the success probabilities of independent Bernoulli random variables, namely as the  indicators
$$I_{j}:=I_{j(m)}~= ~{\bf 1}\,\Big\{ [u_{j-1}, u_j[ {~\rm contains~jump \,times~of~}(C_u)_{\T\le u\le 1}\Big\} $$ for $j=1, 2, \cdots , m.$ The success probability of the $j$th Bernoulli experiment is then given by $p_j=\E(I_j).$ Let us call this interpretation the "Bernoulli model" for increments of the process  $(C_u)$  for the chosen partition of $[\T,1].$ 

To be definite we now confine our interest to equidistant partitions, and in this class to those such that all $p_j <1.$ Let  $$s(m)=\sup _{j \in \{1, 2, \cdots, m\}}
 \{p_{j(m)}\}.$$
 From (11) we obtain $p_j\sim \D_j\eta(u_j)$ and thus, as $\D_j\to 0$, we have $\E(I_j)\to 0$ and also \begin{align}\E(I_{j(m)})\Big/P( \,[u_j, u_{j+1}[ {\rm ~contains~exactly~one~jump~time}) \to 1.\end{align} 

The idea is now the following:  First, if we can interpret any increment $C_{u_{k}}-C_{u_{j}}, j \le k\le m$ as a sum of odds in our Bernoulli models, then the optimal odds-rule for stopping on the last success identifies the optimal rule for stopping on the last sub-interval of the partition containing jump-times of $(C_u).$   Note that for any fixed $m,$ the last Bernoulli success may correspond to more than one point in the last sub-interval containing points (i.e. jump-times of $(C_u)$). Second, in a limiting Bernoulli model defined by letting $m\to \infty,$  the last success corresponds, according to (13), with probability $1$ to the very last jump in $(C_u).$ Hence, provided that the notion of limiting odds is meaningful for the limiting Bernoulli model, the optimal rule in the latter
identifies the optimal rule for stopping on the last jump of $(C_u).$

 We will combine both parts by showing that the continuous-time analogue of  odds in the limiting Bernoulli model {\it is} an intensity measure of a counting process, and we will adapt it to become the process $(C_u)$.

Let $\rho$ be a real-valued non-negative Riemann integrable function $\rho: [0,1] \to \R^+,$ and let
$$\Psi(x, \Delta x):=\int_{x}^{x+\Delta x} \rho(u) 
du.$$ 
\smallskip
We now chose a function $\rho$  in such a way that all $\Psi(u_{j},\Delta_j)$ satisfy the equation
\begin{align} \Psi_j:=\Psi(u_j,\Delta _j)=\frac{p_j}{1-p_j}=r_j,~ j=1,2, \dots, m. \end{align} Note that the existence of such a function $\rho$ is evident for any {\it finite} partition since
the class of Riemann integrable functions contains already infinitely many. If we choose $\rho$ in this class we have $\lim_{\D u \to 0}\Psi(u, \Delta u)/\Delta u$ exists almost everywhere on $[\T,1],$ and this derivative coincides with $\rho(u)$ on $[\T,1].$ 

Now we must check whether such a function $\rho$ exists if we let the mesh size of the partition tend to $0.$ We shall now prove that the function $\rho$ exists and is unique in the limiting Bernoulli model, and that  $\eta$ and $\rho$ coincide almost everywhere on $[\T,1].$ It will thus be justified to call the function $\rho$  the {\it odds-intensity} associated with the (identical) intensity $\eta$ of the process $(C_u)$ on $[\T, 1].$  

 Indeed, recalling $\D_j=1/m$, we will first show that
 
	$${\rm(i)}~~ \rho(u_j)=\lim_{\D_j\to 0}\,\frac{1}{\D_j}r_j=\eta(u_j),~j=1, 2, \cdots$$
	$${\rm(ii)} \lim_{m\to \infty}\,\sum_{j=1}^m\Psi(u_j,\D_j)=\sum_{j=1}^\infty\lim_{m\to \infty}\, \Psi(u_j,\D_j).$$
	
	\bigskip
	
\noindent The limiting equation (i) follows  from the definition of odds in the Bernoulli models, and from (11), since $$\frac{{p_j}}{1-p_j}\frac{1}{\D_j}\sim\frac{1}{\D_j}\frac{\D_j\eta(u_j)}{(1-
{\D_j\eta(u_j)})}=\frac{\eta(u_j)}{1-
{\D_j\eta(u_j)}}\to \eta(u_j)~{\rm as}~\D_j\to 0.$$	
To see (ii), we first recall that for all $j=1,2, \dots, m$ we have $p_j<1$ and thus from (14)
$$p_j\le \Psi_j = p_j/(1-p_j).$$ For fixed $\epsilon$ with $0<\epsilon<1$ we now choose an integer
$m:=m(\epsilon)$ large enough so that  $s(m):=\sup\{p_k: 1\le k \le m\}<\epsilon.$  This is trivially always possible for a finite number $m$ of $p_k$, since, again seen as a function of $\Delta _k,$ we have from (11) that each $p_k\to 0$ as $\Delta_k \to 0+,$ that is, as $m\to \infty. $ Then we obtain
$$p_j\le \Psi_j\le p_j/(1-s(m))\le p_j/(1-\epsilon),$$
or, according to (11) explicitly,
\begin{align}\int_{u_{j-1}}^{u_j} \eta(u)\,du\le \Psi_j\le\frac{1}{1-\epsilon} \int_{u_{j-1}}^{u_j} \eta(u)\,du.\end{align} Since this inequality holds for all $j=1,2, \cdots, m(\epsilon)$, it must hold also for any sum of these terms (column-wise) taken over the same set of indices. In particular this includes 
tail sums beginning at an arbitrary time $x \ge \T.$  Hence, by bounded convergence, (ii) is true.

But then the latter also holds for any random time $x:=\tau \ge \T$, since, by the hypothesis stated in Theorem 2.5, the intensity measure $\eta$ is supposed to be non-random from time $\T$ onwards. Thus for any set of sub-intervals of $[\T,1]$, the limiting odds sum for the limiting Bernoulli model, corresponds to the integral of $\rho$ over the same set of intervals. 
Therefore, in particular, $\rho$ satisfying (14) must satisfy for any $\epsilon$ and equidistant partition with mesh size $\D_j=1/m(\epsilon)$
\begin{align}\E\big(C_1-C_\tau\big)=\int_\tau^1 \eta(u)du\le \int_\tau^1 \rho(u)du\le \frac{1}{1-\epsilon}\E\big(C_1-C_\tau\big).\end{align} 

Since $\epsilon>0$  can be chosen arbitrarily close to $0$ in the inequality (16), it follows from the squeezing theorem 
that the inner integral is bound to coincide with $\E(C_1-C_{\tau}).$ According to (ii), this inner integral is however the limiting tail sum of odds for the limiting Bernoulli model, and (i) implies thus $\rho(u)=\eta(u).$

Finally, letting $\epsilon\to 0+$ in (16) that the inner integral, that is, the limiting tail sum of odds in the limiting Bernoulli model, drops below $1$ if and 
only if
$\E(C_1-C_\tau)$ does so. Hence the proof.\qed

\begin{rem} The preceding criterion is valid independently of whether $\tau$ is a jump-time of $(C_u)$ or not. Indeed, if $\eta(u)>0$ on $[\tau, 1]$ then for all $0<\epsilon<1-\tau$ we have $\E(C_1-C_{\tau+\epsilon})<1.$ Therefore, if $\tau$ happens to be a jump-time of $(C_u)_{\,t\le u \le 1}$ we must also stop on $\tau$.\end{rem}
We are now ready to tackle our main problem.
 
\section{The open question of optimality}
\subsection{Preview and visualisation of our approach}
If the optimal strategy exists, then it must solely be based on all the sequential information we can have, that is, on the information stemming from the history of arrivals (points) and their relative ranks.  
 
Clearly, any strategy is trivially optimal if there are no points so that we can confine our interest to the case $N>0.$ Denote by $N_u$ the number of arrivals up to time $u$. 
If $N>0,$ there is at least one arrival on $[0,1]$, and the first one is a record by definition. 

Due to the i.i.d. structure of points on $[0,1],$ if the decision maker looks back at time $t\in [0,1]$, and if there are preceding arrivals, then he or she knows that their pattern is the outcome of i.i.d. uniformly distributed points on $[0,t].$  The same will hold by looking forward, that is, if there are arrivals then their unordered arrival times are i.i.d. on $[t,1].$ This is true since i.i.d. uniform random variables on a given interval $I$, say, stay i.i.d. conditioned on their location in sub-intervals of $I.$ This is illustrated in the figure below (Fig.1) where arrivals are denoted by *, and where the first * is meant to indicate the arrival time $T_1$.

\bigskip
 $$ |_0............................*...........*...*..\longleftarrow|_t~ .............................~ |_1$$
 $$ |_0............................*...........*...*.....  ~~...|_t\longrightarrow ......................~ |_1$$

 \smallskip
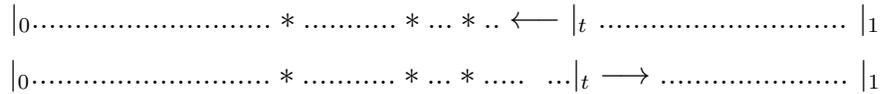
 \centerline{Fig. 1}
 
 \medskip
 \centerline{ Decision-maker's perception}

\bigskip
  
 \medskip
 
 \noindent From the first arrival time $T_1$ onwards ($0<T_1< 1$ a.s.) the decision maker has the information that the counting process  $(N_u)_{u\ge T_1}$ is a process with proportional increments as defined in Subsection 2.6. See Fig. 2.  
Accordingly, given $N_u$, the expected value of the number of points in $[u, u+\Delta u[$ equals $(N_u \Delta u)/u$ almost surely, and it is important to note that no$~o(\Delta u)$ is added here. 
 
 $$~|_0...........................~_{T_1}~...........*...*..(N_u)_{u\ge T_1}......?........?...........|_1$$
  \centerline{Fig. 2}
  
  \medskip
 \centerline{$(N_u)_{u\ge T_1}$ is a 
 proportional-increments process}

 \bigskip \medskip\noindent 
 The
 relevant stochastic process for stopping on rank 1 is then the record process $(R_u)_{u\ge T_1}$ which is a sub-process of $(N_u)_{u\ge T_1}$(see Fig. 3)
 
 \bigskip
  $$~|_0...........................*...........~? \,...\,?~(R_u)_{u\ge T_1} .......?..........?............ |_1$$
  \centerline {Fig. 3}
  
  \medskip
\centerline{$(R_u)_{u\ge T_1}$ is obtained from $
(N_u)_{u\ge T_1}$ by inverse-proportional thinning.}
  
  \bigskip \noindent This thinning is by Rényi's Theorem  such  that if $T_J\ge T_1$ is a jump-time of the process $(N_u)$ then it is retained for the record process $(R_u)_{u\ge {T}_1}$ with probability $1/N_{T_J}$ independently of  retained preceding points. We call this the {\it inverse-proportional thinning} property of R\'enyi's record theorem on the process $(N_u)_{u\ge T_1}.$ 
Note that if we have a predictable non-random intensity measure, the process $(R_u)$ can then play the role of $(C_u)$ in Theorem 2.5.  Stopping online on the desired rank $1$ means stopping online on the very last  record, i.e. on the last jump-time of $(R_u)_{u\ge T_1}.$ 

In the previous paper it was  claimed (see Theorem 3.1) that the $1/e$-strategy is uniquely optimal, but its proof, based on Theorem 2.5, is {\it wrong}. We now point out where exactly the error occurred:

\bigskip
\noindent{\bf ~~~The error in the proof~} 

\smallskip
 We now recapitulate the proof which is correct until equation (22) included :
 
Let ${\cal F}_t$ denote the filtration generated by $\{N_s: 0\le s\le t\},$ and  denote by ${\cal G}_t$ the one generated by both $\{N_s: 0\le s \le t\}$ and $\{R_s: 0\le s \le t\}$ together. Since both fields are clearly   increasing we have ${\cal G}_t\subseteq {\cal G}_u$ for $t\le u \le 1.$ 

Clearly $T_1$ is a $(\G_t)$-measurable stopping time since $(\F_t \subseteq\G_t)$. Given $T_1,$  choose  $t \in [T_1,1]$ and define for fixed $m\in\{2,3, \cdots\}$ and $k=0, 1, 2,  \cdots m-1,$
\begin{align*}  u_k:=u_k(t)=t+\frac{k(1-t)}{m},\\\Delta_k:=\Delta_k(t)=u_{k+1}-u_k=\frac{1-t}{m}.\end{align*} 

\medskip
\noindent It follows that for any ${\cal G}_u$-measurable random variable $X$ and $ 0\le t\le u \le 1,$
\begin{align}  \E(\E(X\big|{\cal G}_u)\,|\,{\cal G}_t)=\E(X\big|{\cal G}_t). \end{align} Let now $X$ denote the number of records in $[t,1]$, that is $X=R_1-R_t.$ 
 From the linearity of the expectation operator we obtain \begin{align}  \E\big(R_1-R_t \big| {\cal G}_t\big)=\E\left(\sum_{k=0}^{m-1} (R_{u_{k+1}}-R_{u_k})\,\Big|\,{\cal G}_t\right)=\sum_{k=0}^{m-1} \,\E\left(R_{u_{k+1}}-R_{u_k}\,\Big|\,{\cal G}_t\right),\end{align} 
and then from (18) used in (17)
\begin{align} 
 \E\big(R_1-R_t \Big| {\cal G}_t\big)=
\sum_{k=0}^{m-1} \,\E\left(\E\left(R_{u_{k+1}}-R_{u_k}\Big|{\cal G}_{u_k}\right)\,\Big|\,{\cal G}_t\right).\end{align}

Let $\lambda(u)$ denote the rate of the point process $(N_t)$ at time $u,$ and  $h(u)$
be the conditional probability of a point appearing at time $u$ being a record. The process $(N_u)$ inherits history-dependence from the p.i.-property so that $\lambda(u)$ is also history-dependent, namely a $\cal F$-predictable intensity process
 for $(N_t)$ relative to the filtration $({\cal F}_t)$. The function $h(u)$ acts like a thinning on the counting process $(N_t),$ retaining only its record-times as events. The resulting record process has an intensity, $\eta$ say, which may depend on both  $\lambda$ and $h$, and which we write formally as
\begin{align}\label{eq12} \eta(u):= g_{\lambda, h}(u)  \,\,,T_1\le u \le1.
\end{align}
Note that this formal definition is a step of caution because $h(u)$ and $\lambda(u)$ are history-dependent random variables, and dependent on each other. Thus we do not assume so far that $\eta(u)= g_{\lambda, h}(u)$ factorises into $\lambda(u)h(u)$ over sub-intervals we will consider.  Of course we know it does so point-wise because $h$ is defined as the conditional probability of a point being retained as a record.

\smallskip Now consider the inner conditional expectation on the r.h.s.\,of (19). Since $(N_u)_{T_1 \le u \le 1}$ is a p.i.-process, and  $N_{u_k}$ is ${\cal G}_{u_k}$-measurable, we have correspondingly $$\E(N_{u_{k+1}}-N_{u_k}\Big|{\cal G}_{u_k})=\Delta_k N_{u_{k}}/u_k ~a.s.,$$ and thus $\lambda(u_k)=N_{u_k}/u_k ~{\rm a.s.}.$  Moreover, if $u_k$ were a jump-time for $(N_u)$ it would be according to R\'enyi's Theorem  a record time with probability $1/N_{u_k}$  which shows that $h$ in (20) is also history-dependent. 

We now show the central fact that the increments of the record process $(R_u)$ on $[u,u+\Delta u[$ given $\G_u$ will never depend on the {\it locations} of jump-times in $[u,u+\Delta u[,$ but only on the number of jumps in there. Indeed, if we denote the $j$th jump-time in $[u_k, u_{k+1}[$ by $A_j:=T_{j+N_{u_k}}$, then
\begin{align*}\E\left (R_{u_{k+1}}-R_{u_k}\Big|N_{u_{k+1}}-N_{u_k}=J;\,A_1, A_2, \cdots , A_J\right)\\=\E\left(\sum_{j=1}^J {\bf 1}\{ A_j {\rm ~is~a ~record~time}\}\right).~~~~~~~~~~ \end{align*}
Since $J\le N_1<\infty $ we can exchange the operators expectation and summation, and then use Rényi's Theorem.
Therefore, by the definition of the $A_j,$ the latter equals
\begin{align} \sum_{j=1}^J P \left( A_j{\rm~is~a~record~time}\,\Big|\, \G_{u_k}\right) = \sum_{j=1}^J\ \frac{1}{N_{u_k}+j} ~{\rm a.s.},\end{align}
which is understood as being zero if $J=0.$
Given $\G _{u_k}$, the value $N_{u_k}$ is a constant, and $J$ is $\F_{u_k}$-predictable. Hence the r.h.s.\,of (20) is $\F_{u_k}$-predictable and {\it does
not} depend on the location of jumps. 

But then, given any interval $[u, u+du]$,  we can imagine these jump times (if any) to be located
where we want them to be within this interval, and we are entitled to think of the first one (if any) as being in $u.$
This implies from (21) that $g_{\lambda, h}$ in (20) must factorise on the sub-interval $[u, u+du]$ into the intensity of $(N_u),$ namely $\lambda(u),$ and the inverse proportional thinning $h(u)=1/N_u.$ Now recall that the p.i.-property of $(N_u)$ for $u\ge T_1$ implies \begin{align}\lambda(u)du:=\E(dN_u|\F_u)=\E(N_{u+du}-N_u\,|\,\F_u)=\frac{N_u}{u}\,du~ {\rm a.s.}. \end{align} Since the inverse-proportional thinning on $(N_u)$ is $(\F_u)$-predictable and $\F_u\subseteq\G_u,$ we have correspondingly
 \medskip

 \centerline{\bf The error was in (23) of [*]. It should read}
\begin{align} \E\left(dR_u\Big|\G_u\right):=\E\left(R_{u+du}-R_u\Big|\G_u\right)=\frac{N_u}{u}\frac {1}{N_u+1}\,du ~{\rm ~a.s.}\end{align} Indeed, with the intensity $N_u/u$ of the arrival process $(N_u)$ given in (22), we have a positive increment of the record process $(R_u)$ if and only if $[u, u+\Delta u]$ contains a jump-time and the latter {\it is}  a record-time which occurs then with probability $du/(N_u+1),$ (and not $du/N_u$.) This means that the stopping time $T_1$ does {\it not} fulfill the condition that $(C_u):=(R_u)$ would have for $u\ge T_1$, independent increments, unless $N_u=\infty.$ Hence the claim is not proved.

 \subsection {Implication of the correction}
 
(I1)~~\noindent If $N_t=\infty$  for some $t>0$ then Theorem 2.5 of [*] can be applied, because then $N_u/(N_u+1)=1$~a.s. Otherwise it cannot be applied since (23) stays history dependent.  However, the case $N_u=\infty$ adds nothing new to what was known before. Indeed, since $N=N_1\ge N_u=\infty$,  we know already from III of Section 2.2 that the $1/e$-strategy is optimal for $N=\infty.$ 

\bigskip\noindent
(I2)~~ Integration of (23) on $[t,1]$ yields
\begin{align*}
\E(R_1-R_t|\G_t) =\int_t^1\E\left(d R_u\Big|\G_u\right)du=\int_t^1 \frac{N_u}{N_u+1} u^{-1}du \le\int_t^1 u^{-1}du = -\log(t),
\end{align*}
where, unless $N_t=\infty$, we have the strict
inequality $\E(R_1-R_t|\G_t)< - \log(t).$ Since $-\log(t)\le 1$ implies $t\le 1/e$ this implies also that, if Theorem 2.5 would apply, the $1/e$-strategy could {\it not} be  optimal in the case $N<\infty ~a.s..$

\bigskip\noindent
(I3)~~\noindent As we have just seen, Theorem 2.5 cannot be 
applied in the case $N<\infty ~a.s.$ Also, we do not know whether the condition $\E(R_1-R_t|\G_t)\le 1$  for a record-time $t$ is at least a {\it necessary} condition for optimal stopping at time $t.$
\smallskip 

Taking I1, I2 and I3 together we conclude that an answer to the open question will need a different approach.  \qed

\section{Optimal strategies without value}
One must be careful in dealing with problems under the hypothesis of no information. Usually, if we speak of a problem of optimal stopping, we think of finding a non-anticipative rule maximizing a pre-determined
objective function, and the solution we find constitutes the value (see e.g. Peskir and Shiryayev (0, Rüschendorf, or Stirzaker (2015).
However, as seen for instance in the l.a.p. of Bruss and Yor (2012), it may occur that a problem of optimal stopping and/or optimal control has no value.  Moreover, as we will show below, it may resist any comparison  of performance versus non-optimal strategies.

\smallskip \noindent The following Lemma illustrates this in a simple form.
\begin{Lem}  In a model for problems of optimal stopping and/or optimal control in a no-information setting, the following features are possible:
\begin{quote}
{\rm(i)} An optimal strategy $\Ess$ solving the defined problem may exist independently of whether one can attribute a value to $\Ess$. 

\medskip
{\rm (ii)} If an optimal strategy $\Ess$ exists, it need not be the limit of $\epsilon$-optimal strategies as $\epsilon\to 0+$.
\end{quote}
\end{Lem}

\begin{rem} In the way Lemma  4.1  is formulated, the statements (i) and (ii) can be proven by examples having properties (i) and (ii).  As said before, the no-information last-arrival problem is such an example. However, the following simple  example suffices to make the point. We keep it in form of a an optimal control problem in order to concentrate on the essence, but by adding costs  for observations we can change the example easily into a stopping problem.\end{rem}

\medskip \noindent
{\bf Proof}

\smallskip
\noindent (i) Let $(I^{(1)}_j)_{j=1, 2, \cdots}$ and $(I^{(2)}_j)_{j=1, 2, \cdots}$ be two sequences of Bernoulli random variables, not necessarily independent of each other, and let $p_j^{(1)}=P(I_j^{(1)}=1)$ and $p_j^{(2)}=P(I_j^{(2)}=1).$ At each time $j$ the decision-maker (he, say) sees both $p_j^{(1)}$ and $p_j^{(2)}$ and decides on which Bernoulli experiment he wants to bet  (see Fig. 4).
If he  bets on Line(1)
he receives the random reward $I_j^{(1)},$ and, alternatively,  if he bets on Line(2), he receives the random reward $I_j^{(2)}$. At time $j$ he sees only the two entries $p_j^{(1)}$ and $p_j^{(2)}$ but none of the future values for $j'>j.$
\begin{align*}{\rm Line~(1):}~~~ p_1^{(1)}~~~p_2^{(1)}~~~p_3^{(1)}~~~\cdots~~~p_j^{(1)}~~~\cdots
~~~p_n^{(1)}\\{\rm Line~(2):}~~~p_1^{(2)}~~~p_2^{(2)}~~~p_3^{(2)}~~~\cdots~~~p_j^{(2)}~~~\cdots~~~p_n^{(2)} \end{align*} 

\medskip
\centerline {Fig. 4}

\bigskip
\medskip\noindent 
Denoting by $ \pi:\N \to\{\text{Line}~(1),\text{Line}~(2)\}$ the decision policy at each step, his objective is to maximise for each $n$ the expected accumulated reward. The optimal strategy, if it exists, is defined by \begin{align}\pi^*=\arg\max_{\pi}\left\{\E\left(\sum_{k=1}^n I_k^{\pi(k)}\right)\right\}.\end{align} But this implies that it does exists: in order to play optimally, it suffices to bet at each step $j$ on $\max\,\{p_j^{(1)}, p_j^{(2)}\}.$ Indeed, this strategy yields at each time $n$ the expected accumulated reward $${\rm M(n)}=\max\,\{p_1^{(1)}, p_1^{(2)}\}+\max\,\{p_1^{(1)}, p_2^{(2)}\}+\cdots+\max\,\{p_n^{(1)}, p_n^{(2)}\}$$ upon which one cannot possibly improve because the maximum of a sum never exceeds the sum of the maxima. And thus we have \begin{align}M(n)=\sum_{k=1}^n \max\left\{p_k^{(1)},p_k^{(2)}\right\}\ge \max_{\pi}\E\left(\sum_{k=1}^n I_k^{\pi(k)}\right).\end{align}
(26) and (27) imply that the optimal strategy ${\Ess}_n$ maximizing the accumulated reward until time $n$ exists, but nevertheless, before time $n,$ no value can be attributed to the optimal ${\Ess}_n$
because $M(n)$ is still unknown. This proves (i).
(We note that if the corresponding values in Line (1) and Line (2) never coincide, all ${\Ess}_n$ are moreover unique. )
 \qed

\medskip
(ii)
To prove (ii), look at the following modification.
Suppose that at some time $t\in \N$ a red light is switched on for Line (2), say, with probability $\delta.$ If the light is switched on, the decision maker is supposed to be no longer entitled to bet on Line (2).  
No information is given how often, or how long, the red light may be switched on, given it is switched on at least once.

It is straightforward to check, similarly as above, that now the unique
optimal strategy is to bet, {\it whenever possible}, on the Line with the entry $\max\{p_j^{(1)}, p_j^{(2)}\}$. If  $\delta=0$ then we are in the case (i). Further we see easily that, if \begin{align}\ell=\lim_{n\to \infty}\,\sum_{j=1}^n \,\left |p_j^{(1)}-p_j^{(2)}\right| < \infty,\end{align}
then, for  any given $\epsilon>0,$ we can always choose $\delta$ sufficiently small so that the optimal strategy  in this setting is $\epsilon$-optimal with respect to $\Ess.$ Indeed, for all $n$ the difference in the accumulated rewards is bounded above by $\delta \ell.$ In this case, the optimal strategy can be seen as the limit of $\epsilon$-optimal strategies.

 If the limit $\ell$ in (28) satisfies $\ell=\infty$, however, then this is not possible. 

In conclusion, we simply do not know whether the existing optimal strategy can be seen as a limit of $\epsilon$-optimal strategies, at least not in this class of $\epsilon$-optimal strategies. 

\smallskip
This does of course  not exclude that one may still be able to find other $\epsilon$-optimal strategies.
However, the point we want to make is that special circumstances in a given problem may {\it naturally} lead us to a certain class of $\epsilon$-optimal strategies with which we would like to study the given problem, because we understand them. Then we would like to  be able to count on some form of {\it closedness} as we know it from other domains of Mathematics.  In {\it Analysis} for instance, we require for good reasons that a function $f: \R^n\to \R$ allows a limit in $x\in \R^n$ if and only if for {\it all} sequences $(x_m)\to x$ we have $f(x_m)\to f(x)$. As we  have just seen, without knowing that $\ell$ defined in (24) satisfies $\ell < \infty,$ we would not know whether all $\epsilon$-optimal strategies would do as $\epsilon \to 0+$. 
 \qed

\subsection {\bf Particularities of the no-information hypothesis}
\smallskip

Lemma 4.1 tells us that we must keep, in more general cases, something important in mind: 
 In a setting of an optimal stopping problem under no-information an optimal strategy need not have a neighbourhood in a classical sense in the set of possible strategies.  An optimal expected payoff
need not be a limit in an analytic sense of the corresponding expected payoffs of seemingly close strategies. 

\smallskip
But then, any argument based on a continuity assumption, or on the existence of a point of indifference for the optimal decision, etc.,  may become questionable. 

\smallskip
This implies that we may,
in certain cases, be able to show the optimality of a certain strategy without
being able to assess  at the same time how a (slightly) sub-optimal strategy, or in fact any other strategy, would compare to the optimal strategy with respect to performance. The non-negligible content of what we point out here is that we have to be careful when speaking about indifference values, limiting performances or even any limit argument in the context of no-information.

\bigskip\bigskip
\centerline{***}
\newpage
\section{References}

~~~~Abdel-Hamid A., Bather J., and Trustrum G. (1982) {\it The secretary problem with
an unknown number of candidates}, J. Appl. Probability, Vol. 19 (3): 619-630.
\medskip

Bruss F.T. (1984)  {\it A unified approach to a class of best choice problems with an
unknown number of options}, Annals of Probability, Vol. 12 (3): 882-889.
\medskip

Bruss F.T. (1987) {\it On an optimal selection problem of Cowan and Zabczyk}, J. Appl. Probability Vol. 24: 918-928.

\medskip

Bruss F.T. (2000) {\it Sum the odds to one and stop}, Annals of Probability, Vol. 28 (3): 1384-1391.
\medskip

Bruss F.T. and Samuels S.M. (1987)), {\it A Unified Approach to a Class of Optimal Selection  Problems with an Unknown Number of Options,} Annals of Probability, Vol. 15:  824-830.

\medskip
Bruss F.T. and Samuels S.M. (1990) {\it Conditions for quasi-stationarity of the Bayes rule in selection problems with an unknown number of rankable options}, Annals of Probability, Vol. 18 (2): 877-886.

\medskip
Bruss F.T. and Rogers L.C.G. (1991) {\it Pascal processes and their characterization},
Stoch. Proc. and Their Applic.,
Vol. 37 (2): 331-338
\medskip

Bruss F.T. and Yor M. (2012) {\it Stochastic processes with proportional increments and the last-arrival problem},  Stoch. Proc. and Their Applic.,
Vol. 122 (9): 3239-3261.

\medskip
Cowan R. and Zabczyk J. (1978) {\it An optimal selection problem associated with the Poisson process}, Theory of Prob. and  Applic., Vol. 23: 548-592.

\medskip Ferguson T.S. (2016) {\it The Sum-the-Odds Theorem with Application to a Stopping game of Sakaguchi},
Mathematica Applicanda, Vol.  44 (1): 45-61.

\medskip
Grau Ribas J.M (2020),{\it An extension of the last-success-problem}, Statistics \& Probability  Letters, Vol. 156, DOI: 10.1016/j.spl.2019.108591 

\medskip
Hadamard  J. (1902) {\it Sur les probl\`emes aux d\'eriv\'ees partielles et leur signification physique.}  Princeton University Bulletin, pp. 49-52.

\medskip
Matsui T. and Ano K. (2016) {\it Lower bounds for Bruss' odds problem with multiple stopping}, ~Math. of Oper. Research, Vol.  41 (2): 700-714.

\medskip
Presman E.L. and Sonin, I.M. (1972) {\it The best choice problem for a random
number of objects}, Theory of Prob. and Applic., Vol. 17 (4): 657-668.

\medskip
R\'enyi A. (1962) {\it Th\'eorie des \'el\'ements saillants d'une suite d'observations}, Annales
scientifiques de l'Universit\'e de Clermont-Ferrand 2, S\'erie Math\'ematiques, 8 (2): 7-13.

\medskip

R\"uschendorf  L. (2016)
{\it Approximative solutions of optimal stopping and selection problems}, Mathematica Applicanda, Vol. 44 (1): 17-44.

\medskip
Samuels, S.M. (1985)  Math Reviews: MR0744243 (85m:62182).
\medskip

Stewart T.J. (1981) {\it The secretary problem with an unknown number of options.}
Operations Research, Vol. 29 (1): 130-145.
\medskip

Stirzaker D. (2015) {\it The Cambridge Dictionary of  Probability and Its Application}, Cambridge University Press, ~ISBN 978-1-107-07516-0.

\bigskip\bigskip

\noindent{\bf Author's address}
\medskip:\\ F.\,Thomas Bruss, \\Universit\'e Libre de Bruxelles, CP 210, \\B-1050 Brussels, Belgium\\ (tbruss@ulb.ac.be)

\end{document}